\documentclass{amsart}
\usepackage{times}
\usepackage[T1]{fontenc}
\usepackage[english]{babel}
\usepackage{latexsym,amssymb}

\author{davide bondoni}
\title{Schr\"oder e la teoria dei gruppi}

\begin{document}
\maketitle

\section{Introduzione}
Nel 1874 Ernst Schr\"oder (1841--1902), allora insegnante al {\it Realgymnasium} in Baden-Baden, pubblica come "supplemento" un piccolo libricino dal titolo {\it Sugli elementi formali di un'algebra assoluta}.\footnote{\cite{algebra}.} In questo lavoro, si suole riconoscere la nascita dell'{\it algebra universale}. Fra l'altro, ha suscitato molto interesse, la netta distinzione e chiarificazione tra componente {\it sintattica} e componente {\it semantica} in una teoria formale, tanto che recentemente Volker Peckhaus (Padeborn) ha parlato al proposito di {\it nascita della teoria dei modelli}.\\
\indent Due aspetti colpiscono il lettore:
\begin{enumerate}
\item in quel periodo Schr\"oder non aveva nessuna cognizione logica; tutto il suo background era puramente ed esclusivamente matematico. Schr\"oder si muover\`a ad occuparsi di logica solo nel 1876--1877;
\item la letteratura su Schr\"oder, sempre molto filosofica e logica, non ha prestato molto interesse agli aspetti pi\`u tecnici e matematici di {\it Sugli elementi formali}.
\end{enumerate} 
Il nostro interesse si concentrer\`a sul punto (2), rimandando il lettore interessato ad una visione pi\`u logica alla letteratura sull'argomento.

\section{Algoritmi}
Schr\"oder introduce un insieme $I$ di elementi qualsiasi su cui impone della {\it struttura}. Ovvero, su $I$ definisce una serie di operazioni ed individua degli elementi {\it privilegiati}. Queste strutture egli le chiama con un nome alquanto sfortunato, {\it algoritmi}\footnote{Ad essere sinceri, Schr\"oder usa l'espressione 'algoritmo' in maniera ambigua; una volta indicando un {\it insieme struttuato} in un certo modo, una volta indicando le {\it operazioni} definite su quell'insieme.} e gli elementi privilegiati come {\it moduli} (elementi neutri, nel nostro vocabolario). In realt\`a, con {\it algoritmo}, Schr\"oder non intende una regola computazionale, ma un insieme su cui sono state definite alcune operazioni. Aprendo il testo a p.~14, troviamo elencati questi {\it algoritmi} $\forall a,b \in I$:
\begin{align}
ab = ba,\qquad a:b=\frac{a}{b}\qquad (C_1)\notag\\
a:b = b:a, \qquad \frac{b}{a} = ba,\qquad (C_2)\notag\\
\frac{a}{b}=\frac{b}{a}, \qquad ab=b:a, \qquad (C_3)\notag\\
ab=a:b,\qquad a:b=\frac{b}{a}, \qquad \frac{b}{a} = ab, \qquad (C_0)\notag
\end{align}
Alla semplice giustapposizione di elementi $a,b \in I$ corrisponde il prodotto tra $a$ e $b$, talora indicato con $a.b$. La prima equazione di $C_1$ ci dice che due qualsiasi elementi di $I$ commutano tra loro rispetto al prodotto ($ab=ba$). La seconda equazione ci dice che la divisione {\it sinistra} (indicata con $:$) coincide con quella {\it destra} (indicata dalla frazione $\frac{a}{b}$). In $I$ strutturato da $C_1$ c'\`e quindi {\it una} sola forma di divisione.\\
\indent La prima equazione di $C_2$ ci dice che due qualsiasi elementi di $I$ commutano tra loro rispetto alla divisione sinistra, la seconda che la divisione destra coincide con la moltiplicazione. Si noti che in $C_2$ sia il prodotto che la divisione destra possono {\it non} essere commutativi.\\
\indent In $C_3$ gli elementi di $I$ commutano tra loro rispetto alla divisione destra, per il resto, il prodotto coincide con la divisione sinistra. Riguardo alla commutativit\`a di queste due operazioni, vale quanto osservato sopra a proposito di $C_2$. In breve, ognuno dei primi tre algoritmi si caratterizza per almeno una operazione commutativa.\\
\indent Con $C_0$, invece, non \`e commutativa {\it nessuna} operazione e tutte e tre le operazioni coincidono (banalmente, per transitivit\`a, si ha che da $ab = a:b$ e $a:b = \frac{b}{a}$ segue che $ab = \frac{b}{a} = a:b$).\\
\indent A questo gruppo di equazioni, Schr\"oder aggiunge dei {\it moduli}, ovvero, $\forall a,b \in I$:
\begin{align}
a.a = b.b,\qquad (M_x)\notag\\
a:a = b:b,\qquad (M_:)\notag\\
\frac{a}{b}= \frac{b}{a},\qquad (M_/)\notag
\end{align}
Questi moduli corrispondono a ci\`o che noi chiameremmo oggi {\it elementi neutri}, rispettivamente del prodotto, della divisione {\it sinistra} e di quella {\it destra}. Definendo gli algoritmi $C_0,...,C_3$ sul nostro insieme di partenza $I$ che dotiamo anche dei tre moduli $M_x,M_:$ e $M_/$, otteniamo quattro tipi di strutture:
\begin{align}
\mathcal{G}_1=<I,.,M_x>\notag\\
\mathcal{G}_2=<I,:,M_:>\notag\\
\mathcal{G}_3=<I,/,M_/>\notag
\end{align}
Come \`e facile vedere, si tratta di tre {\it gruppi} commutativi, o abeliani, e associativi. A fare eccezione \`e la quarta struttura che consiste in un {\it gruppoide} e le cui operazioni {\it non} sono commutative:
\begin{equation}\notag
\mathcal{G}_0 = <I,.,:,/>
\end{equation}

\section{Interpretazione}
Come detto, Schr\"oder non parla di {\it gruppi}, ma di {\it algoritmi}, anche se \`e chiaro che quando si tratta di considerare gli algoritmi come strutture, essi coincidono con la nostra nozione di {\it gruppo}.\\
\indent La nozione di {\it gruppo} era stata introdotta per la prima volta da \'Evariste Galois nella sua famosa {\it prima memoria},\footnote{{\it M\'emoire sur les conditions de r\'esolubilit\'e des \'equations par radicaux} (1831--1832). Ora in \cite[106--144]{neumann}.} per indicare (con una certa ambiguit\`a) un insieme di {\it permutazioni} o di {\it sostituzioni}, in funzione della dimostrazione che un'equazione di quinto grado non pu\`o essere risolta mediante radicali. Un gruppo, per Galois, quindi, era un qualcosa formulato per uno scopo ben preciso. Si potrebbe osservare che secondo questa prospettiva, qualsiasi altro concetto matematico atto allo scopo, per Galois sarebbe andato ugualmente bene. A Galois non interessava affatto indagare il concetto di gruppo in s\'e e per s\'e. Tanto \`e vero, che confonde un gruppo di permutazioni con uno di sostituzioni.\footnote{\cite[p.~20 e segg.]{neumann}.} Inoltre, la sua concezione puramente combinatoria o computazionale della matematica gli impedivano una tale concezione {\it strutturale} degli oggetti matematici. Per usare le stesse parole di Galois:
\begin{quotation}
[La vita della scienza] \`e bruta e assomiglia a quella dei minerali che crescono per giustapposizione\ldots Gli analisti non deducono, \textbf{combinano}, \textbf{compongono}: per quanto immateriale possa essere, l'analisi non \`e in nostro potere pi\`u di qualsiasi altra [scienza]; bisogna spiarla, sondarla, sollecitarla.\footnote{Citato in \cite[p.~262]{goldstein}. Il grassetto \`e mio.}
\end{quotation}
\indent Schr\"oder conosceva le opere di Galois? Mi sento di rispondere con franchezza: no, nella maniera pi\`u assoluta. Ma conosceva i testi di Serret e Jordan che contribuirono non poco a divulgare il nuovo concetto di {\it gruppo}. Tuttavia, non \`e possibile evincere neppure da questi testi un'indagine meta-teorica del concetto di gruppo. Bisogner\`a aspettare Felix Klein che introdurr\`a quella visione strutturale della teoria dei gruppi a noi cos\`i familiare e celebrata da van der Waerden nel suo capolavoro {\it Algebra Moderna}.\\
\indent Schr\"oder non ritorner\`a pi\`u sull'argomento, ma il modo in cui tratta i suoi algoritmi \`e squisitamente strutturale. Egli confronta i vari algoritmi in base alla loro potenza, cio\`e al numero di conseguenze che possono essere tratte da loro. Cos\`i facendo, sposta l'attenzione dal concetto di algoritmo, nel senso comune della parola, come qualcosa di utile ad uno scopo determinato, per intrapprendere un'indagine sul significato stesso del concetto di algoritmo/gruppo. Per esempio, l'algebra assoluta menzionata nel titolo, risulta essere un gruppo commutativo abeliano che estende l'usuale algebra (anche perch\'e le operazioni sono puramente formali).\\
\indent Di fatto, Schr\"oder dedica solo poche pagine all'argomento, forse perch\'e incapace di realizzare la portata di ci\`o che stava facendo, ma \`e un peccato che queste pagine siano finora sfuggite agli storici impegnati nella teoria dei gruppi o in quella di Galois. Si tratta infatti della prima indagine meta-teorica (strutturale) del concetto di {\it gruppo}, condotta con mezzi combinatori.

\nocite{jordan}
\nocite{serret}
\nocite{waerden}

\bibliographystyle{amsalpha}
\bibliography{pristem}

\end{document}